\documentclass[11pt]{amsart}
\usepackage[a4paper,margin=27mm]{geometry}
\usepackage[T1]{fontenc}
\usepackage{lmodern}
\usepackage{amsmath,amssymb,amsthm,mathtools}
\usepackage{microtype}
\usepackage[colorlinks=true,linkcolor=blue,citecolor=blue,urlcolor=blue]{hyperref}
\usepackage{enumitem}
\setlist[enumerate]{itemsep=2pt,topsep=4pt}
\newtheorem{theorem}{Theorem}[section]
\newtheorem{proposition}[theorem]{Proposition}
\newtheorem{lemma}[theorem]{Lemma}
\newtheorem{corollary}[theorem]{Corollary}
\theoremstyle{remark}
\newtheorem{remark}[theorem]{Remark}
\newcommand{\R}{\mathbb R}
\newcommand{\C}{\mathbb C}
\newcommand{\Sp}{\operatorname{Sp}}
\newcommand{\U}{\operatorname{U}}
\newcommand{\SU}{\mathrm{SU}}
\newcommand{\SO}{\mathrm{SO}}
\newcommand{\diag}{\operatorname{diag}}
\newcommand{\Tr}{\operatorname{Tr}}
\newcommand{\Ad}{\operatorname{Ad}}

\newcommand{\ip}[2]{\langle #1,#2\rangle}
\title[Thompson's Diagonal Theorems]{Geometric Proofs of Thompson's Diagonal Theorems}

\author{Yongzhi Luan}
\address{School of Mathematics, Shandong University, Jinan, Shandong 250100, China}
\email{luanyongzhi@email.sdu.edu.cn}
\author{Tin-Yau Tam}
\address{Department of Mathematics and Statistics, University of Nevada, Reno, NV 89557, USA}
\email{ttam@unr.edu}

\date{\today}


\begin{document}

\begin{abstract}
We give a conceptual proof of Thompson's characterization of the diagonals
of complex symmetric matrices with prescribed singular values. The
nonnegative diagonal-modulus set is the moment polytope for the restriction
of a compact symplectic coadjoint orbit to $\Sp(1)^n$; this identification
follows from the real convexity theorem of O'Shea and Sjamaar. The necessary
inequalities follow from symplectic linear algebra and the decomposability
of bivectors in dimension three. A relative Schubert-degree argument, using
the subgroup moment-cone criterion of Berenstein and Sjamaar, proves
sufficiency. The exceptional inequality
has two complementary explanations: decomposability in dimension three
proves its necessity, and the exceptional Schubert-degree identity $3=1+2$ accounts for
its place in the sufficient inequality list. A companion section gives a short
geometric proof of the exceptional inequality in Thompson's earlier theorem
for general complex matrices; sufficiency there follows from Kostant's
type-$D_n$ convexity theorem. For both matrix classes, fixing the sum
of the diagonal moduli gives the convex hull of the permutations of a
single explicit vector. We use these sections to determine the vertices
and to compare the two sets of admissible moduli. Finally, for the diagonal-torus problem under special-unitary equivalence, we prove that Thompson's 1994 conjectured necessity condition holds in orders at most three and give a counterexample in order four.
\end{abstract}
\maketitle

\section{The theorem and the geometric question}

Let $s_1\ge\cdots\ge s_n\ge0$, write $S=\diag(s_1,\ldots,s_n)$, and put
\[
 \sigma=\sum_{i=1}^n s_i,\qquad A_k=\sum_{i=1}^k s_i,\qquad A_0=0.
\]
Takagi factorization identifies the complex symmetric matrices with these
singular values as $USU^T$, $U\in\U(n)$. Define
\[
 D(s)=\{\diag(USU^T):U\in\U(n)\},\qquad
 M(s)=\{(|b_{11}|,\ldots,|b_{nn}|):B=USU^T\}.
\]
The coordinates in $M(s)$ are not ordered. Diagonal unitary congruence
allows their phases to be prescribed independently, so determining $M(s)$
determines $D(s)$.

\begin{theorem}[Thompson \cite{Thompson1979}]\label{thm:Thompson}
Let $d_1,\ldots,d_n\in\C$, and arrange $r_i=|d_i|$ so that
$r_1\ge\cdots\ge r_n$. There is a complex symmetric matrix with diagonal
$d$ and singular values $s$ if and only if
\begin{align}
 \sum_{i=1}^k r_i&\le A_k &&(1\le k\le n),\tag{T1}\label{eq:T1}\\
 \sum_{i=1}^{k-1}r_i-\sum_{i=k}^n r_i
 &\le\sigma-2s_k &&(1\le k\le n),\tag{T2}\label{eq:T2}\\
 \sum_{i=1}^{n-3}r_i-\sum_{i=n-2}^n r_i
 &\le\sigma-2(s_{n-1}+s_n) &&(n\ge3).\tag{T3}\label{eq:T3}
\end{align}
The last condition is omitted for $n=1,2$, and an empty sum is zero.
\end{theorem}

In his note added in proof, Thompson \cite[p.~105]{Thompson1979}
suggested that Lie theory should explain the result without an elaborate
analysis of cases. Tam \cite[Section~6]{Tam1999} obtained a Lie-theoretic
description of the real-part projection of $D(s)$ and noted at that time
the absence of a conceptual proof of the full theorem. The real-part
projection is governed by weak majorization: if $x\in\R^n$, then
$x\in\operatorname{Re}D(s)$ precisely when
$\sum_{j=1}^k|x|_j^\downarrow\le A_k$ for all $k$. Here the arrow
denotes decreasing rearrangement. This projection does not retain
\eqref{eq:T2}--\eqref{eq:T3}.

We instead keep the absolute values of the diagonal entries. The Takagi
orbit is realized as the fixed locus of an antisymplectic involution on
a compact $\Sp(n)$-orbit. Projection to $\Sp(1)^n$ records a vector in
each three-dimensional block, whose length on the fixed locus is
$|b_{jj}|$. The real convexity theorem of O'Shea and Sjamaar identifies
these lengths with a convex moment polytope. Trace estimates prove
\eqref{eq:T1}--\eqref{eq:T3}; the relative Schubert-degree criterion then
shows that these inequalities suffice.

The two Thompson theorems have a parallel symmetric-space origin.
The isotropy representation of $\Sp(2n,\R)/\U(n)$ identifies its
noncompact tangent space with complex symmetric matrices, on which
$\U(n)$ acts by $B\mapsto UBU^T$. For general complex matrices the
corresponding space is $\SU(n,n)/\mathrm{S}(\U(n)\times\U(n))$, with action
$A\mapsto UAV^*$. Kostant's linear projection in the latter setting
records real diagonal components, whereas compressing the Hermitian
dilation of $A$ to coordinate two-planes records $|a_{ii}|$ as a
positive eigenvalue. Our compact $\Sp(n)$-orbit uses analogous
two-plane coordinates for the symmetric case. That coadjoint orbit
is an auxiliary space associated with the involution, rather than
the symmetric space $\Sp(2n,\R)/\U(n)$ itself. The real-part
projection associated with $\SU(n,n)$ (restricted type $C_n$) is
distinct from the type-$D_n$ projection on real matrices used in
Section~\ref{sec:general-complex}.

We state the required real convexity theorem and verify its hypotheses
in Section~\ref{sec:orbit}. For the symmetric theorem, the other
substantial external input is the subgroup moment-cone
criterion of Berenstein--Sjamaar \cite{BerensteinSjamaar2000}, used
in Section~\ref{sec:schubert-sufficiency}. Its nonvanishing condition
is used only through the necessary equality of Schubert degrees; no Schubert
coefficients are computed. For a complementary geometric invariant theory
approach to the generalized eigenvalue problem, see Ressayre
\cite{Ressayre2010}. The matrix, exterior-algebra, and
majorization arguments are proved below.

In Section~\ref{sec:general-complex} we return to Thompson's 1977 theorem
for general complex matrices. Its exceptional inequality follows from a
one-dimensional intersection in a Hermitian dilation. Together with the
type-$D_n$ Kostant description proved in \cite[Section~2]{Tam1999}, this
gives a short conceptual proof of that theorem as well.
Section~\ref{sec:mass-sections} then uses the two theorems to describe
the sections of fixed total diagonal modulus and to measure the effect
of imposing symmetry. Finally, Section~\ref{sec:su-torus}
returns to Thompson's subsequent diagonal-torus problem under
special-unitary equivalence. We settle his conjectured necessity
condition for matrices of order three and exhibit its failure in order four.

For later use, let $P(s)\subset\R_{\ge0}^n$ be the permutation-invariant
polytope defined by the inequalities in Theorem~\ref{thm:Thompson}.
Equivalently, impose \eqref{eq:T1} for every subset of $k$ coordinates,
\eqref{eq:T2} for every subset of $k-1$ plus signs, and \eqref{eq:T3}
for every subset of $n-3$ plus signs. This is a finite linear system in
the joint variables $(s,r)$ on the chamber $s_1\ge\cdots\ge s_n\ge0$.

\section{A compact symplectic orbit containing the Takagi orbit}
\label{sec:orbit}

We first recall Takagi factorization in the form needed here. For a
complex symmetric matrix $B$, the antilinear map $Tv=B\overline v$
satisfies $T^2=B\overline B=BB^*$ and
$\|Tv\|^2=\langle v,T^2v\rangle$. It preserves each eigenspace of
$T^2$. On an eigenspace with eigenvalue $s^2>0$, the map $s^{-1}T$
is an antiunitary involution and therefore has an orthonormal basis
of fixed vectors: its fixed real subspace spans the eigenspace over
$\C$, and real Gram--Schmidt applies there. On the zero eigenspace
$T=0$. Combining these bases gives a unitary matrix $U$ such that
$B\overline U=US$, or $B=USU^T$, with the singular values on the
diagonal of $S$.

Set $V=\C^n\oplus\C^n$ and
\[
 J_0=\begin{pmatrix}0&I_n\\-I_n&0\end{pmatrix},\qquad
 G=\Sp(n)=\{g\in\U(2n):g^T J_0g=J_0\}.
\]
We use the Hermitian realization $\mathfrak p_G=i\mathfrak{sp}(n)$
of the dual of the compact Lie algebra, with pairing
$\ip{A}{C}=\frac12\Tr(AC)$. Its elements are
\begin{equation}\label{eq:Hermitian-model}
 A=\begin{pmatrix}H&B\\\overline B&-\overline H\end{pmatrix},
 \qquad H^*=H,\quad B^T=B.
\end{equation}
Indeed these are precisely the Hermitian matrices satisfying
$A^TJ_0+J_0A=0$. Their eigenvalues occur in opposite pairs. Let
$\mathcal O_s$ denote the $G$-orbit of
\[
 A_S=\begin{pmatrix}0&S\\S&0\end{pmatrix}.
\]
It is a compact connected coadjoint orbit, with eigenvalues
$s_1,\ldots,s_n,-s_n,\ldots,-s_1$ and its usual
Kirillov--Kostant--Souriau symplectic form. More explicitly, for
skew-Hermitian $X,Y\in\mathfrak{sp}(n)$ we use the convention
\[
 \Omega_A([X,A],[Y,A])=-\frac{\mathrm i}{2}\Tr(A[X,Y]).
\]
Then $\ell_C(A)=\frac12\Tr(CA)$ is the Hamiltonian of the
infinitesimal action generated by $\mathrm iC$, for Hermitian
$C\in\mathfrak p_G$. Restricting this identity to a subgroup shows
that its moment map is orthogonal projection onto its Hermitian
Lie algebra model.

Write $f_i=e_{n+i}$ in the standard basis of $V$. Let
$E_i=\operatorname{span}_{\C}\{e_i,f_i\}$ be the $i$th coordinate
symplectic two-plane and set
\[
 H_0=\Sp(1)^n\subset G.
\]
To avoid confusing this subgroup with the matrix $H$ in
\eqref{eq:Hermitian-model}, we denote it throughout by $H_0$.
Its moment map is the orthogonal projection
\begin{equation}\label{eq:block-moment}
 \mu(A)=\bigl(P_{E_i}A|_{E_i}\bigr)_{i=1}^n.
\end{equation}
The $i$th block is
\[
 \begin{pmatrix}h_{ii}&b_{ii}\\\overline{b_{ii}}&-h_{ii}\end{pmatrix},
\]
whose eigenvalues are $\pm\rho_i$, where
$\rho_i=(h_{ii}^2+|b_{ii}|^2)^{1/2}$. Define
\[
 \Delta(s)=\{(\rho_1,\ldots,\rho_n):A\in\mathcal O_s\}.
\]
The chamber of $H_0$ is $\R_{\ge0}^n$, without an ordering condition.

We record the special case of real convexity used in the proof.
\begin{theorem}[O'Shea--Sjamaar {\cite[Theorem~3.1]{OSheaSjamaar2000}}]
\label{thm:real-convexity}
Let $L$ be a compact connected Lie group, with an involution
$\vartheta$, acting on a compact connected Hamiltonian manifold $X$
with equivariant moment map $\mu$. Suppose that an antisymplectic involution
$\tau$ of $X$ has a nonempty fixed locus and satisfies
\[
 \tau(lx)=\vartheta(l)\tau(x),\qquad
 \mu\circ\tau=-\vartheta\circ\mu.
\]
Identify the dual Lie algebra using an $\Ad(L)$- and
$\vartheta$-invariant inner product. Suppose a Cartan subalgebra
$\mathfrak t$ is contained in the negative eigenspace of $\vartheta$, and let
$\mathfrak t_+$ be a positive Weyl chamber in $\mathfrak t$. Then
\[
 \mu(X^\tau)\cap\mathfrak t_+
       =\mu(X)\cap\mathfrak t_+,
\]
and this set is a compact convex polytope.
\end{theorem}

\begin{proposition}\label{prop:real-convexity}
For every $s_1\ge\cdots\ge s_n\ge0$,
\[
                         M(s)=\Delta(s).
\]
In particular, $M(s)$ is a compact convex polytope.
\end{proposition}
\begin{proof}
Put $K_0=\diag(iI_n,-iI_n)$ and define
$\vartheta(g)=K_0gK_0^{-1}$. This is an involutive automorphism of $G$,
and its fixed group is the embedded $\U(n)$ consisting of
$\diag(U,\overline U)$. On \eqref{eq:Hermitian-model}, $\vartheta$
keeps $H$ and changes the sign of $B$. The map
\[
                  \tau(A)=-\vartheta(A)
\]
preserves $\mathcal O_s$, fixes $A_S$, and satisfies
$\tau(gA_Sg^*)=\vartheta(g)A_S\vartheta(g)^*$.
The displayed formula for $\Omega$ and invariance of trace give
$\tau^*\Omega=-\Omega$. Thus $\tau$ is antisymplectic, and it obeys
$\tau(gAg^*)=\vartheta(g)\tau(A)\vartheta(g)^*$.
Its fixed locus is
\[
 \mathcal O_s^\tau
 =\left\{\begin{pmatrix}0&B\\\overline B&0\end{pmatrix}:
                       B=USU^T,\ U\in\U(n)\right\}.
\]
Here the square of the displayed matrix has diagonal blocks
$BB^*$ and $B^*B$, so its positive eigenvalue parameters are exactly
the singular values of $B$; Takagi factorization proves the equality.

Apply Theorem~\ref{thm:real-convexity} to the $H_0$-action on
$\mathcal O_s$. Both the group and the orbit are compact and connected,
and $A_S$ is a fixed point of $\tau$. Compactness makes the moment map
proper, as required in the general O'Shea--Sjamaar theorem.
The subgroup $H_0$ is $\vartheta$-stable and
$\mu\tau=-\vartheta\mu$. In each $\Sp(1)$ factor the off-diagonal
real line
\[
                \begin{pmatrix}0&r\\r&0\end{pmatrix},\quad r\in\R,
\]
is both a maximal abelian subspace of the negative eigenspace and a
Cartan subalgebra. Thus $\mathfrak a=\mathfrak t$ for $H_0$.
The positive chamber corresponds to $r_i\ge0$. The chamber image of
$\mathcal O_s^\tau$ is $M(s)$: diagonal unitary congruence makes every
$b_{ii}$ nonnegative, independently. The full chamber image is
$\Delta(s)$. Real convexity proves the assertion.
\end{proof}

After this identification, we may conjugate by $H_0$ so that the
compression to $E_i$ is $\diag(r_i,-r_i)$. This changes neither its
chamber coordinate nor the ambient orbit. This diagonal Cartan
convention will be used in the rest of the proof.

\section{Necessity from symplectic linear algebra}
\label{nec:section}

This section supplies the necessary inequalities directly. In particular,
it does not appeal to Thompson's theorem or to a Schubert-calculus
criterion. Let $V=\mathbb C^{2n}=E_1\oplus\cdots\oplus E_n$ be the
orthogonal decomposition into coordinate symplectic two-planes. Write
$\omega(v,w)=v^TJ_0w=\sum_i\omega_i(v,w)$ for the standard complex
symplectic form. A complex Lagrangian is an $n$-dimensional subspace
on which $\omega$ vanishes. Let $Jv=J_0\overline v$ be the compatible antiunitary
operator, with $J^2=-1$. A Hermitian matrix $A$ in the symplectic Lie algebra model
satisfies $AJ=-JA$ and has eigenvalues $\pm s_1,\ldots,\pm s_n$.
For now suppose $s_n>0$. The relation $AJ=-JA$ gives
$JE_-=E_+$ for the negative and positive spectral spaces. Since these
spaces are Hermitian orthogonal and $\dim E_\pm=n$, they are
Lagrangian. More generally, for an $n$-plane $L$, isotropy is equivalent
to $JL=L^\perp$. A sum of one line from each $E_i$ is Lagrangian.

Conjugate by $H_0=\operatorname{Sp}(1)^n$ so that the compression of $A$
to $E_i$ has eigenvalues $r_i,-r_i$, with $r_i\geq0$. Let $\ell_i^+$
be a positive eigenline of that compression; when $r_i=0$, choose any
line. For a subset $Q\subseteq\{1,\ldots,n\}$, consider the product
Lagrangians
\begin{equation}\label{nec:product}
 L=\bigoplus_{i\notin Q}\ell_i^+\ \oplus\!\bigoplus_{i\in Q}\ell_i,
 \qquad \ell_i\subset E_i\text{ any complex line}.
\end{equation}
Every such $L$ satisfies
\begin{equation}\label{nec:lowertrace}
 \operatorname{Tr}(AP_L)
 \geq\sum_{i\notin Q}r_i-\sum_{i\in Q}r_i,
\end{equation}
where $P_L$ is the Hermitian orthogonal projection. Indeed, its trace
is the sum of the Rayleigh quotients on the mutually orthogonal
coordinate lines, each varying between $-r_i$ and $r_i$.

We use the following elementary trace bounds. If a Hermitian operator
$T$ has eigenvalues $\lambda_1\ge\cdots\ge\lambda_N$ and $P$ is
an orthogonal projection of rank $k$, then
\begin{equation}\label{eq:trace-variational}
 \sum_{j=N-k+1}^N\lambda_j\le\Tr(TP)
                   \le\sum_{j=1}^k\lambda_j.
\end{equation}
Indeed, in an orthonormal eigenbasis $v_j$,
$\Tr(TP)=\sum_j\lambda_j t_j$, where
$t_j=\|Pv_j\|^2\in[0,1]$ and $\sum_jt_j=k$. Moving weight toward
the largest, respectively smallest, eigenvalues proves the bounds.
These are the Ky Fan trace inequalities.

For any set $I$ of
$k$ coordinates, the positive compression eigenlines in those planes
form an orthogonal $k$-plane. Its trace is $\sum_{i\in I}r_i$, so the
Hermitian Ky Fan variational principle gives
\begin{equation}\label{nec:T1}
                 \sum_{i\in I}r_i\le A_k.
\end{equation}
Here and below the trace on a subspace means $\Tr(AP_L)$; the subspace
need not be invariant under $A$.

\begin{lemma}\label{nec:tracebound}
If a Lagrangian $L$ contains a subspace $R\subseteq E_-$, then
\begin{equation}\label{nec:traceupper}
 \operatorname{Tr}(AP_L)
 \leq\sigma-2\operatorname{Tr}(A_-P_R),
\end{equation}
where $A=A_+-A_-$ is the positive/negative spectral decomposition
and $\sigma=\sum_i s_i$.
\end{lemma}
\begin{proof}
Since $L$ is Lagrangian, $JL=L^\perp$, so $JR\subseteq E_+$ is
orthogonal to $L$. Consequently
\[
 \operatorname{Tr}(A_+P_L)
 \leq\operatorname{Tr}(A_+)-\operatorname{Tr}(A_+P_{JR}),
 \qquad
 \operatorname{Tr}(A_-P_L)\geq\operatorname{Tr}(A_-P_R).
\]
The identity $AJ=-JA$ gives
$\operatorname{Tr}(A_+P_{JR})=\operatorname{Tr}(A_-P_R)$.
Subtracting proves \eqref{nec:traceupper}.
\end{proof}

\begin{proposition}\label{nec:T2}
For every $k=1,\ldots,n$ and every subset $Q$ of size $n-k+1$,
\begin{equation}\label{nec:ineqT2}
 \sum_{i\notin Q}r_i-\sum_{i\in Q}r_i\leq\sigma-2s_k.
\end{equation}
\end{proposition}
\begin{proof}
Let $F_k\subset E_-$ be spanned by eigenvectors with eigenvalues
$-s_1,\ldots,-s_k$, and put
\[
 V_0=\bigoplus_{i\notin Q}\ell_i^+\ \oplus\!\bigoplus_{i\in Q}E_i.
\]
Since $\dim V_0=n+|Q|$ and $\dim F_k=k=n-|Q|+1$, the intersection
$F_k\cap V_0$ contains a nonzero vector $v$. For each $i\in Q$ choose
$\ell_i$ to contain the $E_i$-component of $v$, choosing an arbitrary
line if that component vanishes. The resulting Lagrangian
\eqref{nec:product} contains $v$. With $R=\mathbb Cv$ and $v$
normalized, $\langle v,A_-v\rangle\geq s_k$. Combine
\eqref{nec:lowertrace} and \eqref{nec:traceupper}.
\end{proof}

\begin{lemma}\label{nec:threeplanes}
Suppose $Q$ has three elements. There is a product Lagrangian of the
form \eqref{nec:product} whose intersection with $E_-$ has dimension
at least two.
\end{lemma}
\begin{proof}
The subspace $V_0$ in the preceding proof has dimension $n+3$,
so $E_-\cap V_0$ has dimension at least three. Choose a three-dimensional
subspace $W$ of this intersection. Pull the three coordinate
symplectic forms $\omega_i$, $i\in Q$, back to $W$ by the coordinate
projections. Their sum is zero: $W\subset E_-$ is isotropic, while
all components outside $Q$ lie in the fixed one-dimensional lines
$\ell_i^+$ and make no contribution to $\omega|_W$.

Thus these three alternating forms impose at most two independent
linear conditions on the three-dimensional space $\Lambda^2W$.
Their common kernel contains a nonzero bivector $\xi$. Every such
bivector is decomposable. Indeed, after permuting a basis,
$\xi=a u_1\wedge u_2+b u_1\wedge u_3+c u_2\wedge u_3$ with
$a\ne0$, and
\[
 \xi=\left(u_1-\frac ca u_3\right)\wedge(a u_2+b u_3).
\]
Write $\xi=u\wedge v$ with $u,v$ independent, and set
$R=\operatorname{span}\{u,v\}\subset W$.
For each $i\in Q$, $\omega_i(\pi_i u,\pi_i v)=0$, where $\pi_i$ is
the coordinate projection to $E_i$. Since $E_i$ is a symplectic
two-plane, this says that the projection of $R$ to $E_i$ has dimension
at most one. Choose $\ell_i$ containing that projection. The product
Lagrangian \eqref{nec:product} then contains $R\subset E_-$.
\end{proof}

\begin{proposition}\label{nec:T3}
For $n\geq3$ and every subset $Q$ of size three,
\begin{equation}\label{nec:ineqT3}
 \sum_{i\notin Q}r_i-\sum_{i\in Q}r_i
 \leq\sigma-2(s_{n-1}+s_n).
\end{equation}
\end{proposition}
\begin{proof}
Choose $L$ and a two-dimensional $R\subseteq L\cap E_-$ from
Lemma~\ref{nec:threeplanes}. The restriction of $A_-$ to $E_-$ has
eigenvalues $s_1,\ldots,s_n$, so the lower Ky Fan variational bound
gives $\operatorname{Tr}(A_-P_R)\geq s_{n-1}+s_n$.
Apply \eqref{nec:lowertrace} and Lemma~\ref{nec:tracebound}.
\end{proof}

Together with the ordinary Ky Fan upper bounds, these are all of
Thompson's necessary inequalities. For a zero singular value, write
$A=gA_Sg^*$, replace $s$ by $s+\varepsilon(1,\ldots,1)$, and keep
$g$ fixed. The moment block norms and all right sides converge as
$\varepsilon\downarrow0$, proving necessity for every $s$.
The exceptional inequality has an elementary geometric
source: in a three-dimensional space, every bivector is decomposable,
and the three coordinate alternating forms satisfy one linear relation.

\begin{remark}[Hollow matrices in orders three and four]\label{rem:sel-hollow}
There is also a matrix explanation of the distinction at three; see
Thompson \cite[Lemma~9 and the subsequent remark]{Thompson1979}.
If $C=C^T$ is hollow of order three, choose real phases satisfying
$\theta_i+\theta_j=-\arg c_{ij}$ for $i<j$, assigning an arbitrary
argument to a zero entry. These three equations are nonsingular.
The diagonal unitary congruence with entries $e^{i\theta_i}$ makes
$C$ real symmetric and still hollow. Its eigenvalues sum to zero, so
the largest of their absolute values is the sum of the other two:
$t_1=t_2+t_3$ for its singular values. Conversely, when this equality
holds, the real matrix $\diag(t_1,-t_2,-t_3)$ has trace zero and is
orthogonally similar to a hollow matrix. Indeed a traceless real
symmetric form has a unit vector of zero Rayleigh quotient; repeat
on its traceless orthogonal compression.

In order four the condition is instead $t_1\le t_2+t_3+t_4$.
Necessity is \eqref{eq:T2} at $k=1$. For sufficiency, the intervals
$[t_1-t_2,t_1+t_2]$ and $[t_3-t_4,t_3+t_4]$ meet. A common length
gives complex $z_j$ with $|z_j|=t_j$ and $z_1+z_2=-(z_3+z_4)$.
Put $Q=\tfrac12\left(\begin{smallmatrix}1&1\\1&-1\end{smallmatrix}\right)^{\otimes2}$.
Then $Q$ is real orthogonal, and every diagonal entry of
$Q\diag(z_1,z_2,z_3,z_4)Q^T$ is $\tfrac14\sum_j z_j=0$.
This is the phase construction of Wootters
\cite[eqs.~(23)--(24)]{Wootters1998}. Thus $(1,1,1,0)$ admits a
hollow symmetric matrix although $(1,1,1)$ does not.
\end{remark}

\section{Sufficiency from relative Schubert degrees}
\label{sec:schubert-sufficiency}

In this section assume $s_1>\cdots>s_n>0$ and $n\ge2$.
The geometric identification of Section~\ref{sec:orbit} gives
$M(s)=\Delta(s)$, the $\operatorname{Sp}(1)^n$ chamber image of the
$\operatorname{Sp}(n)$ coadjoint orbit. We use the scalar form of the
subgroup moment-cone criterion of Berenstein and Sjamaar
\cite[Theorem~3.3.1]{BerensteinSjamaar2000}. The criterion is a
necessary and sufficient family of inequalities. We will need only
the equality of cohomological degrees required by its nonvanishing
condition; the Schubert coefficients themselves need not be computed.

\begin{lemma}\label{schubert:candidates}
If $r\in\R_{\ge0}^n$ does not belong to $\Delta(s)$, there are vectors
$c,\beta\in\{0,1,-1\}^n$ with the same number $m>0$ of nonzero
entries such that
\begin{equation}\label{schubert:violated}
 c\cdot r>\sum_{k=1}^n\beta_k s_k,
\end{equation}
and the following condition holds. If $q$ is the number of $-1$
entries of $c$, $a_{>k}$ and $z_{>k}$ count respectively the
nonzero and zero entries of $\beta$ after $k$, and $z_{<k}$ counts
its zero entries before $k$, then
\begin{equation}\label{schubert:degree}
 q=J(\beta):=
 \sum_{\beta_k\ne0}z_{<k}
 +\sum_{\beta_k=-1}(a_{>k}+1+2z_{>k}).
\end{equation}
\end{lemma}
\begin{proof}
Use the common maximal torus of
$H_0=\operatorname{Sp}(1)^n$ and $G=\operatorname{Sp}(n)$.
The positive chamber of $H_0$ is the nonnegative orthant, while
the positive chamber of $G$ consists of nonnegative coordinates
in decreasing order. The latter chambers subdivide the former
according to coordinate order. Their rays are the vectors
$\chi_I=\sum_{i\in I}e_i$ for nonempty subsets $I$.
The centralizers of a ray with $|I|=m$ give the partial flag
varieties $(\mathbb{CP}^1)^m$ for $H_0$ and the isotropic
$m$-plane Grassmannian for $G$.

The longest Weyl elements of both groups act by $-1$. Hence the
slice of the moment cone in the cited criterion is the chamber
image $\Delta(s)$ of the orbit with parameter $s$. Its scalar
inequality for a ray $\chi_I$ and a pair of Weyl elements is
\[
 \ip{\widetilde w^{-1}r+w^{-1}s}{\chi_I}\ge0.
\]
The $H_0$ Weyl element independently changes coordinate signs.
Specifically, set $c=-\widetilde w\chi_I$ and $\beta=w\chi_I$,
using the Euclidean identification of the common torus and its dual.
Then this inequality has the form $c\cdot r\le\sum_k\beta_k s_k$,
where the $m$ nonzero entries of $c$ record its support and
$\beta_k$ records whether $s_k$, $-s_k$, or neither is assigned
to the active symplectic coordinates.

To specify the Schubert indices, put $\zeta_m=e_1+\cdots+e_m$
and choose a permutation $v_I$ such that $v_I\zeta_m=\chi_I$.
The partial-flag map is
\[
 \phi_I:(\mathbb{CP}^1)^I\longrightarrow\operatorname{IG}(m,2n),
 \qquad (L_i)_{i\in I}\longmapsto\bigoplus_{i\in I}L_i .
\]
Let $a$ and $u$ be the shortest elements in the cosets
$wv_I\operatorname{Stab}_{W_G}(\zeta_m)$ and
$\widetilde w_0\widetilde w\operatorname{Stab}_{W_{H_0}}(\chi_I)$,
respectively, where $\widetilde w_0$ is the longest element for $H_0$.
Then $a\zeta_m=\beta$ and $u\chi_I=c$. In the conventions of
\cite[Section~3.3]{BerensteinSjamaar2000} the nonzero pairing is
$\langle\phi_I^*\sigma_a,\widetilde c_u\rangle\ne0$.
The cocycle $\sigma_a$ is evaluation-dual to the same-index Schubert
cycle and has degree $2\ell(a)$; the source homology class
$\widetilde c_u$ has degree $2\ell(u)=2q$, since $u$ flips exactly
the $q$ negative active coordinates of $c$. Hence $\ell(a)=q$.
For completeness, count its inversions using the roots
$2e_i$, $e_i-e_j$, and $e_i+e_j$.
An active entry preceded by a zero contributes one; each
negative active entry contributes its self root, one root
for each later active entry, and two for each later zero.
The count is $\ell(a)=J(\beta)$ in \eqref{schubert:degree}, so
$J(\beta)=q$. This is a necessary condition
for nonvanishing. Below we impose every inequality satisfying it,
including those whose Schubert coefficients vanish. Berenstein--Sjamaar
therefore supplies a violated inequality of the stated form whenever
$r\notin\Delta(s)$.
\end{proof}

We next show that all the candidate inequalities in this lemma
follow from Thompson's three families. Put $p=m-q$ and $z=n-m$.
For any normal $c$, the Ky Fan inequalities give
\begin{equation}\label{schubert:kyfan}
 c\cdot r\le\sum_{c_i=1}r_i\le A_p.
\end{equation}

\begin{lemma}\label{schubert:mixedzeros}
If $q>0$, $z>0$, and $J(\beta)=q$, then
$\sum_k\beta_k s_k\ge A_p$.
\end{lemma}
\begin{proof}
Set $J=J(\beta)$ and $p_0=m-J=p$. We claim that every prefix
of $\beta$ satisfies
\begin{equation}\label{schubert:prefix}
 \sum_{k=1}^j\beta_k\ge\min(j,p_0),\qquad 1\le j\le n.
\end{equation}
A prefix of only positive entries satisfies this. Otherwise,
let $A'$ be the number of active entries after the prefix and
$L$ the number of negative entries in it. Its sum is
$m-A'-2L$, so it suffices to show $J\ge A'+2L$.
If the prefix contains a zero, its root pairs with the $A'$
later active positions give at least $A'$ inversions. Its
pairs with the $L$ negative positions in the prefix give
at least $L$ more, and the negative self roots give another
$L$. If the prefix has a negative but no zero, fix a zero
after it. Its pairs with the $L$ earlier negatives give
$2L$ inversions, their self roots give $L$, and their pairs
with later active positions give $LA'$. Thus
$J\ge3L+LA'\ge A'+2L$. These counted roots are distinct.
In either case the prefix sum is at least $p_0$. Since a prefix
that is not entirely positive has sum at most $j-1$, this forces
$p_0<j$, and hence proves \eqref{schubert:prefix}.

Abel summation and $s_1>\cdots>s_n>0$ give
\[
 \sum_k\beta_k s_k-A_{p_0}
 =\sum_{j=1}^{n-1}(s_j-s_{j+1})
  \sum_{k=1}^j(\beta_k-\mathbf1_{k\le p_0})
  +s_n\sum_{k=1}^n(\beta_k-\mathbf1_{k\le p_0})\ge0.
\]
\end{proof}

\begin{lemma}\label{schubert:distinctparts}
Let $\mathcal A$ be a set of at least two distinct positive
integers, let $q=\sum_{a\in\mathcal A}a$, and suppose
$0\le t_1\le\cdots\le t_q$. If $\mathcal A\ne\{1,2\}$, then
\begin{equation}\label{schubert:tailmajorization}
 \sum_{j=1}^q t_j\ge2\sum_{a\in\mathcal A}t_a.
\end{equation}
\end{lemma}
\begin{proof}
Write $\mathcal A=\{a_1<\cdots<a_l\}$ and set $t_0=0$.
The difference between the two sides of
\eqref{schubert:tailmajorization} is
\[
 \sum_{j=1}^q(t_j-t_{j-1})
 \bigl(q-j+1-2\#\{a\in\mathcal A:a\ge j\}\bigr).
\]
The coefficients in parentheses decrease between successive elements
of $\mathcal A$ and jump upward immediately after each such element;
past $a_l$ they are positive. Thus it suffices to check
nonnegativity at $j=a_h$. If $l\ge3$, then
\[
 q-a_h+1=\sum_{i\ne h}a_i+1
 \ge\frac{l(l+1)}2-h+1\ge2(l-h+1).
\]
For $l=2$, the two required bounds are $a_2+1\ge4$ and
$a_1+1\ge2$; only $\{a_1,a_2\}=\{1,2\}$ fails the first.
\end{proof}

\begin{theorem}\label{schubert:sufficiency}
If $s_1>\cdots>s_n>0$, every nonnegative vector satisfying
\eqref{eq:T1}--\eqref{eq:T3} belongs to $M(s)$, with
\eqref{eq:T3} omitted for $n=2$.
\end{theorem}
\begin{proof}
We show that every inequality permitted by
Lemma~\ref{schubert:candidates} follows from the stated
families. If $q=0$, the degree formula forces all active
assignments to be positive and to precede all zero entries.
Then $\beta=(1,\ldots,1,0,\ldots,0)$ and its inequality
is a Ky Fan inequality of type \eqref{eq:T1}.
If $q>0$ and $z>0$, Lemma~\ref{schubert:mixedzeros} and
\eqref{schubert:kyfan} imply its inequality.

Suppose $z=0$. Put $N=\{k:\beta_k=-1\}$. If $N$ is empty,
the candidate follows from the total-sum case of
\eqref{eq:T1}. Otherwise the degree condition reduces to
\begin{equation}\label{schubert:partition}
 q=\sum_{k\in N}(n-k+1),
 \qquad b=\sum_k\beta_k s_k
       =\sigma-2\sum_{k\in N}s_k.
\end{equation}
For $N=\{k\}$, this says $q=n-k+1$ and $p=k-1$,
so the candidate is exactly an inequality of type
\eqref{eq:T2}.

If $|N|\ge2$, take $\mathcal A=\{n-k+1:k\in N\}$ and
$t_j=s_{n-j+1}$ for $1\le j\le q$. Unless
$N=\{n-1,n\}$, Lemma~\ref{schubert:distinctparts} gives
$\sum_{j=1}^q t_j\ge2\sum_{k\in N}s_k$ and therefore
$b\ge A_{n-q}=A_p$. The candidate follows from
\eqref{schubert:kyfan}. The remaining assignment has
$N=\{n-1,n\}$ and $q=2+1=3$; its inequality is
exactly \eqref{eq:T3} for a set of three negative coordinates.

All the inequalities required by the Berenstein--Sjamaar
criterion hold, so $r\in\Delta(s)=M(s)$.
\end{proof}

The exceptional three-minus inequality arises from the unique
multiple-negative assignment not absorbed by Ky Fan. The
decomposable-bivector argument of
Lemma~\ref{nec:threeplanes} supplies its direct geometric
necessity. Permutation symmetry puts the three minus signs
on the smallest diagonal moduli when the inequalities are
written with ordered coordinates.

\section{Repeated singular values, zero singular values, and phases}

We finish the proof without assuming regularity of $s$. For $n=1$,
\eqref{eq:T1} and \eqref{eq:T2} say $r_1=s_1$, and the assertion is
immediate. Assume $n\ge2$.

The preceding sections prove necessity for every $s$, and sufficiency
when $s_1>\cdots>s_n>0$. The passage to a boundary spectrum must
approximate the pair $(s,r)$; perturbing $s$ alone can lose feasibility.
Take $s^\circ=(n,n-1,\ldots,1)$ and $r^\circ=s^\circ$, realized by
the diagonal matrix. By necessity, $(s^\circ,r^\circ)$ satisfies the
joint linear system defining $P$. If $r\in P(s)$, then
\[
 s^{(\varepsilon)}=(1-\varepsilon)s+\varepsilon s^\circ,
 \qquad
 r^{(\varepsilon)}=(1-\varepsilon)r+\varepsilon r^\circ
 \quad(0<\varepsilon<1)
\]
also satisfy that system. The spectrum $s^{(\varepsilon)}$ is strictly
decreasing and positive. Regular sufficiency therefore supplies
$U_\varepsilon\in\U(n)$ with
\[
 \left|\diag\bigl(U_\varepsilon
       \diag(s^{(\varepsilon)})U_\varepsilon^T\bigr)\right|
       =r^{(\varepsilon)}.
\]
Compactness of $\U(n)$ gives a convergent subsequence as
$\varepsilon\downarrow0$, and its limit realizes $r$ at spectrum $s$.
Thus $P(s)\subset M(s)$ for every $s$, while necessity gives the
reverse inclusion.

Finally, if a symmetric matrix $B$ has $|b_{ii}|=|d_i|$, choose
unit complex numbers $z_i$ with $z_i^2b_{ii}=d_i$; if $b_{ii}=0$,
then $d_i=0$ and $z_i$ is arbitrary. Congruence by
$\diag(z_1,\ldots,z_n)$ preserves the singular values and gives
the required diagonal. A permutation congruence restores any original
coordinate order. This completes the proof of
Theorem~\ref{thm:Thompson}.

\section{The nonsymmetric theorem and a one-free-plane argument}
\label{sec:general-complex}

The same spectral two-plane viewpoint gives a short proof of the exceptional
inequality in Thompson's earlier theorem for general complex matrices.
This is exactly the $k=n$ instance of \eqref{eq:T2}; the argument below
establishes it without assuming symmetry.
The sufficiency argument then comes from the type-$D_n$ Kostant polytope
in \cite[Section~2]{Tam1999}.

\begin{theorem}[Thompson \cite{Thompson1977}]
\label{thm:general-complex}
Let $d_1,\ldots,d_n\in\C$, and arrange their moduli as
$r_1\ge\cdots\ge r_n$. There is a complex $n\times n$ matrix with
diagonal $d$ and singular values $s_1\ge\cdots\ge s_n\ge0$ if and only if
\begin{align}
 \sum_{i=1}^k r_i&\le\sum_{i=1}^k s_i
       &&(1\le k\le n),\label{eq:general-kyfan}\\
 \sum_{i=1}^{n-1}r_i-r_n
       &\le\sum_{i=1}^{n-1}s_i-s_n.\label{eq:general-exceptional}
\end{align}
\end{theorem}

\begin{proof}
For $n=1$, the two inequalities say $r_1=s_1$. Assume $n\ge2$.
To prove necessity, let $A$ be a matrix with the prescribed data.
After a simultaneous permutation of rows and columns, assume
$|a_{ii}|=r_i$. Form its Hermitian dilation
\[
 X=\begin{pmatrix}0&A\\A^*&0\end{pmatrix}
 \quad\text{on}\quad
 \C^{2n}=E_1\oplus\cdots\oplus E_n,
 \qquad E_i=\operatorname{span}_{\C}\{e_i,e_{n+i}\}.
\]
The eigenvalues of $X$ are
$s_1,\ldots,s_n,-s_n,\ldots,-s_1$, and its compression to
$E_i$ has eigenvalues $r_i,-r_i$. Choose a $+r_i$ eigenline
$\ell_i^+\subset E_i$, with any line permitted when $r_i=0$.
For any subset $I$ of $k$ coordinates, the sum of the lines
$\ell_i^+$, $i\in I$, is a $k$-plane of trace $\sum_{i\in I}r_i$.
The upper Ky Fan bound gives \eqref{eq:general-kyfan}.

For \eqref{eq:general-exceptional}, leave the last coordinate plane free:
\[
 V_0=\ell_1^+\oplus\cdots\oplus\ell_{n-1}^+\oplus E_n,
 \qquad \dim V_0=n+1.
\]
Choose an $n$-dimensional spectral subspace $F_-$ of $X$ with
eigenvalues $-s_1,\ldots,-s_n$. If some $s_i$ vanish, choose the
needed vectors from the zero eigenspace. Since
$\dim V_0+\dim F_-=2n+1$, take a unit vector
$v\in V_0\cap F_-$. Let $\ell_n\subset E_n$ contain the
$E_n$-component of $v$, choosing any line if that component is zero.
Then $L=\ell_1^+\oplus\cdots\oplus\ell_{n-1}^+\oplus\ell_n$
is an $n$-plane containing $v$. The compression bounds on its
coordinate lines give
\[
 \sum_{i=1}^{n-1}r_i-r_n\le\operatorname{Tr}(XP_L).
\]
On the other hand, $\ip{v}{Xv}\le-s_n$ and the upper Ky Fan bound
on the $(n-1)$-plane $L\ominus\C v$ give
\[
 \operatorname{Tr}(XP_L)
 =\ip{v}{Xv}+\operatorname{Tr}(XP_{L\ominus\C v})
 \le-s_n+\sum_{i=1}^{n-1}s_i.
\]
This proves the exceptional inequality.

For sufficiency, let $W_{D_n}$ be the group of signed permutations
with an even number of sign changes. Here the real symmetric space
is $\SO_0(n,n)/(\SO(n)\times \SO(n))$. Its tangent space at the base
point consists of blocks
$\left(\begin{smallmatrix}0&B\\ B^T&0\end{smallmatrix}\right)$,
$B\in M_n(\R)$; the compact group acts by $B\mapsto UBV^T$.
Projection onto the diagonal Cartan subspace records $\diag B$.
The type-$D_n$ application of
Kostant's linear convexity theorem in \cite[Section~2]{Tam1999} gives
\begin{equation}\label{eq:real-kostant-diagonals}
 \left\{\diag(USV^T):U,V\in \SO(n)\right\}
       =\operatorname{conv}(W_{D_n}s),
 \qquad S=\diag(s_1,\ldots,s_n).
\end{equation}
For dominant $x$ and $s$, membership of $x$ in this Weyl polytope
is equivalent to the inequalities against the fundamental weights.
Indeed, if $u^+$ is the dominant Weyl translate of a normal $u$,
the Weyl rearrangement inequality gives
$\ip{u}{x}\le\ip{u^+}{x}$, while the support function of the
polytope at $u$ is $\ip{u^+}{s}$. Thus a separating normal may be
taken dominant, and the dominant cone is generated by the fundamental
weights. They are
\begin{align*}
 \varpi_k&=e_1+\cdots+e_k&& (1\le k\le n-2),\\
 \varpi_{n-1}&=\tfrac12(e_1+\cdots+e_{n-1}-e_n),&
 \varpi_n&=\tfrac12(e_1+\cdots+e_n).
\end{align*}
For the nonnegative decreasing vectors $r$ and $s$, their inequalities
are \eqref{eq:general-kyfan} for $k\le n-2$, together with
\[
 \sum_{i=1}^{n-1}r_i-r_n\le\sum_{i=1}^{n-1}s_i-s_n,
 \qquad
 \sum_{i=1}^{n}r_i\le\sum_{i=1}^{n}s_i.
\]
The $k=n-1$ Ky Fan inequality is the average of the last two.
The same two terminal-weight formulas apply when $n=2$, where
$D_2\cong A_1\times A_1$. Thus the hypotheses put
$r\in\operatorname{conv}(W_{D_n}s)$, so
\eqref{eq:real-kostant-diagonals} supplies a real matrix with diagonal
$r$ and singular values $s$. Left multiplication by a diagonal unitary
matrix prescribes the complex phases of $d$, and a permutation of rows
and columns restores the original coordinate order.
\end{proof}

The free two-plane makes $\dim V_0=n+1$, forcing it to meet the
spectral $n$-plane for $-s_1,\ldots,-s_n$. This explains why one negative diagonal
modulus appears in Thompson's exceptional inequality.
The exceptional inequality is also the terminal type-$D_n$ weight
inequality with one minus sign. For symmetric matrices, the Lagrangian
relation yields the remaining inequalities in \eqref{eq:T2}; the stronger
two-dimensional intersection in the symplectic three-plane argument
yields \eqref{eq:T3}.

\section{Fixed total diagonal modulus and the effect of symmetry}
\label{sec:mass-sections}

How does symmetry change the set of admissible diagonals at the same
singular spectrum? Theorems~\ref{thm:Thompson} and~\ref{thm:general-complex}
let us compare the two sets one hyperplane at a time. This comparison
also arises for reciprocal and unrestricted scattering matrices
\cite{GuoFan2022}. Write $M_{\rm sym}(s)=M(s)$ and let $M_{\rm gen}(s)$
be its unrestricted counterpart. Keep the coordinates labelled, and put
$m=\sum_i r_i$. For $0\le p<n$, let $t=(s_{p+1},\ldots,s_n)$, $q=n-p$,
and $\delta=t_1-\sum_{i=2}^q t_i$. Define $\eta_p$ by the relevant
column of the table:
\[
\begin{array}{c|cc}
 q &\text{symmetric}&\text{unrestricted}\\ \hline
 1&t_1&t_1\\
 2&t_1-t_2&0\\
 3&|\delta|&0\\
 q\ge4&\max(0,\delta)&0 .
\end{array}
\]
Write $R_p(r)$ for the sum of the $p$ largest entries of $r$, with
$R_0=0$. In either class the inequalities already proved are equivalent to
\begin{equation}\label{eq:sel-uniform}
 r\ge0,\qquad R_k(r)\le A_k\ (1\le k\le n),\qquad
 2R_p(r)-m\le A_p-\eta_p\ (0\le p<n).
\end{equation}
Indeed $2R_p-m\le A_p$ already follows from $R_p\le m,A_p$.
In the symmetric case \eqref{eq:T2} with $k=p+1$ has right
side $A_p-\delta$; when $q=3$, \eqref{eq:T3} gives $A_p+\delta$
as well. Taking the smallest of the available bounds gives the table.
In the unrestricted case only $q=1$ requires the improvement from
\eqref{eq:general-exceptional}.

\begin{theorem}[Sections of fixed total modulus]\label{thm:sel-slice}
In either class the attainable totals $m$ form $[\eta_0,\sigma]$.
For each such $m$ the section $\sum_i r_i=m$ is the convex hull of
the coordinate permutations of one decreasing vector $z(m)$.
Its sum is $m$, and its proper prefix sums are
\begin{equation}\label{eq:sel-prefix}
 H_p(m)=\min\{A_p,m,(m+A_p-\eta_p)/2\}\qquad(1\le p<n).
\end{equation}
\end{theorem}
\begin{proof}
For $n=1$ the claim is immediate. For $n\ge2$ put
$L_p=A_p-\eta_p$ and $U_p=A_p+\eta_p$. We first check
$U_p\le L_{p+1}$ for $0\le p<n-1$, or equivalently
$\eta_p+\eta_{p+1}\le s_{p+1}$. In the general case this follows
from the table. In the symmetric case, for a tail of length two it is
equality, and for a tail $(a,b,c)$ it is
$|a-b-c|+(b-c)\le a$ since $a\ge b\ge c$.
For a longer tail with first entry $a$, if $\eta_p=0$ use
$\eta_{p+1}\le s_{p+2}\le a$. Otherwise $a-\eta_p$ is the sum of
the remaining entries, which bounds $\eta_{p+1}$.

The total of any admissible vector lies between $U_0=\eta_0$ and
$U_{n-1}=\sigma$ by \eqref{eq:sel-uniform}. In this interval define
$z(m)$ by the following two formulas; $0^j$ denotes $j$ zeros.
On $U_p\le m\le L_{p+1}$, $0\le p<n-1$, put
\[
 z(m)=(s_1,\ldots,s_p,\ m-A_p,\ 0^{n-p-1}).
\]
On $L_p\le m\le U_p$, $1\le p<n$, put
\[
 b=(m-L_p)/2,\qquad
 z(m)=(s_1,\ldots,s_{p-1},\ s_p-\eta_p+b,\ b,\ 0^{n-p-1}).
\]
These intervals cover $[\eta_0,\sigma]$ and the formulas agree at
shared endpoints. In the first formula $0\le m-A_p\le s_{p+1}$;
in the second $0\le b\le\eta_p\le s_{p+1}\le s_p$. Thus both
vectors are decreasing, nonnegative, and of sum $m$, even at repeated
or zero singular values.

For each $j$, $H_j(m)$ equals $m$ up to $L_j$, $(m+L_j)/2$ between
$L_j$ and $U_j$, and $A_j$ thereafter. The ordering of the intervals
therefore gives the prefix sums \eqref{eq:sel-prefix} in both formulas.
Thus a single vector attains the sharp upper bound for every
prefix at once.
So $z(m)$ satisfies \eqref{eq:sel-uniform} and is realizable by the
two diagonal theorems, in their respective classes.

A vector of the same sum satisfies \eqref{eq:sel-uniform} precisely
when it is nonnegative and its decreasing prefix sums are at most
those of $z(m)$. This is precisely the permutation hull of $z(m)$:
the prefix bounds hold on that hull; conversely, rearranging any
linear functional in decreasing order and using summation by parts
bounds its value by its maximum at a permutation of $z(m)$.
Separation gives the reverse inclusion.
\end{proof}

In particular, $\eta_p$ is the least possible total diagonal modulus
for the tail $(s_{p+1},\ldots,s_n)$, by applying the theorem in
dimension $n-p$. The absolute value in the three-entry tail records
which of \eqref{eq:T2} and \eqref{eq:T3} is tighter; the two bounds
agree exactly when $s_{n-2}=s_{n-1}+s_n$.

\begin{corollary}[Vertices]\label{cor:sel-vertices}
For $n\ge2$ the vertices of either modulus polytope are the coordinate
permutations of the following vectors, with duplicates removed:
\begin{align*}
 u_p&=(s_1,\ldots,s_p,\eta_p,0^{n-p-1})       &&(0\le p<n),\\
 v_p&=(s_1,\ldots,s_{p-1},s_p-\eta_p,0^{n-p}) &&(1\le p<n).
\end{align*}
The symmetric list has at most $2n-2$ permutation orbits. In the
unrestricted list the candidates are the
truncations $(s_1,\ldots,s_p,0^{n-p})$ for $0\le p\le n-2$,
the full vector $s$, and $(s_1,\ldots,s_{n-2},s_{n-1}-s_n,0)$.
\end{corollary}
\begin{proof}
The displayed points are the endpoints of the affine pieces of $z(m)$.
Every permuted piece joins the same permutation of its two endpoints.
The section formula and convexity thus give their convex hull.
Conversely, at $u_p$ the first $p$ coordinate-prefix equalities
$\sum_{i\le j}r_i=A_j$ and the last $n-p-1$ zero coordinates fix all
but one entry. The signed equality
$\sum_{i\le p}r_i-\sum_{i>p}r_i=A_p-\eta_p$ fixes the remaining
entry. These are supporting equalities by \eqref{eq:sel-uniform}.
At $v_p$ use the first $p-1$ prefix equalities, the last $n-p$
zero coordinates, and the same signed equality. Hence every point
listed is extreme. In the symmetric case $v_{n-1}=u_{n-2}$;
in the unrestricted case $\eta_p=0$ for $p\le n-2$.
\end{proof}

\begin{corollary}[Where symmetry can be detected]\label{cor:sel-distance}
Let $n\ge2$, $\mathcal S=M_{\rm sym}(s)$, $\mathcal G=M_{\rm gen}(s)$,
and now take $\eta_p$ from the symmetric column. The projections of
$\mathcal S$ and $\mathcal G$ onto any proper subset of coordinates
coincide. For $x\in\mathcal G$ their difference is measured by
\[
 d(x)=\max\left(0,\ \max_{0\le p\le n-2}
     \bigl\{2R_p(x)-\textstyle\sum_i x_i-A_p+\eta_p\bigr\}\right).
\]
The $\ell^1$ distance from $x$ to $\mathcal S$ is $d(x)$.
When $d(x)>0$, any one prescribed coordinate can be changed by
$d(x)$, keeping the others fixed, to obtain a point of $\mathcal S$.
The $\ell^1$ Hausdorff distance between the two polytopes is
$\max_{0\le p\le n-2}\eta_p$.
\end{corollary}
\begin{proof}
Project onto $k<n$ coordinates. The Ky Fan upper bounds hold for
the projection of $\mathcal G$. For any linear functional of these
coordinates, sort its coefficients and let $j$ be the number of
positive ones. Summation by parts and nonnegativity show that its
maximum under those upper bounds is attained at
$(s_1,\ldots,s_j,0^{k-j})$. This is a projection of a permutation
of the symmetric $u_j$: put its remaining entry $\eta_j$ in an
omitted coordinate. Thus the two compact convex projections have
equal support functions.

For $I\subset\{1,\ldots,n\}$, write $c_I=+1$ on $I$ and $-1$
elsewhere. By \eqref{eq:sel-uniform}, the only cuts on $\mathcal G$
needed to give $\mathcal S$ are
$c_I\cdot x\le A_p-\eta_p$ for $p=|I|\le n-2$; the case
$p=n-1$ is already common to both. The greatest positive excess over
these cuts, or zero if none is violated, is $d(x)$. Two distinct cuts
cannot both be violated: their normals have opposite signs at some
coordinate $j$, whereas the projection result gives a point of
$\mathcal S$ agreeing with $x$ away from $j$. No change of that one
coordinate can repair two violations requiring opposite directions.

If $x\notin\mathcal S$, let $c_I$ be the unique violated normal. For
each $j$ the projection result now gives $y=x-t(c_I)_j e_j\in\mathcal S$
with $t\ge d(x)$. At $t=d(x)$ the violated cut becomes tight; all
other cuts hold on the segment $[x,y]\subset\mathcal G$. This gives
the required correction. Conversely $c_I\cdot(x-z)\ge d(x)$ for all
$z\in\mathcal S$, so $\|x-z\|_1\ge d(x)$. Finally, on $\mathcal G$
one has $c_I\cdot x\le A_p$, and equality is attained at the
unrestricted truncation of length $p$ supported in $I$. Taking the
maximum over $p$ gives the Hausdorff distance, since
$\mathcal S\subseteq\mathcal G$.
\end{proof}

For example, at $s=(1,1,1,1)$ the vector $(1,0,0,0)$ occurs for a
general matrix. Its distance to the symmetric set is one: changing
the first entry to zero or any one of the other three entries to one
gives a symmetric diagonal-modulus vector of the same spectrum.

\section{The diagonal torus under special-unitary equivalence}
\label{sec:su-torus}

The preceding sections concern the diagonals attainable under unitary
congruence or unrestricted unitary equivalence. Thompson
\cite{Thompson1994} subsequently investigated a different question: when
does the special-unitary orbit of a given matrix realize \emph{every}
diagonal vector with prescribed absolute values? This distinction is
significant because individual diagonal phases cannot necessarily be
assigned using a determinant-one left multiplier.

For $A\in\C_{n\times n}$ and $r=(r_1,\ldots,r_n)\in\R_{\ge0}^n$, put
\[
 \mathcal T(r)=\{(r_1\lambda_1,\ldots,r_n\lambda_n):
                      |\lambda_1|=\cdots=|\lambda_n|=1\}.
\]
We say that the $\SU(n)\times \SU(n)$ orbit of $A$ \emph{covers}
$\mathcal T(r)$ if, for each $d\in\mathcal T(r)$, there are
$U,V\in \SU(n)$ with $\diag(UAV)=d$. If the singular values of $A$
are $s_1\ge\cdots\ge s_n\ge0$ and $r_1\ge\cdots\ge r_n$, Thompson
proved that
\begin{equation}\label{eq:SU-sufficient}
 \sum_{j=1}^k r_j\le \sum_{j=1}^k s_j\quad(1\le k<n),\qquad
 \sum_{j=1}^n r_j\le \sum_{j=1}^{n-1}s_j-s_n
\end{equation}
is sufficient for full torus coverage. He conjectured that the final
inequality is also necessary \cite[p.~969]{Thompson1994}.
The initial partial-sum inequalities are necessary by the ordinary
Ky Fan bound. We show that the conjectured final inequality is indeed
necessary in order three but fails in order four.

\subsection{Necessity in order three}

The key observation is an elementary inequality for three complex
numbers whose sum is real and whose product is negative real.

\begin{lemma}\label{lem:SU-three-complex}
Let $z_1,z_2,z_3\ne0$ be complex numbers such that
$z_1+z_2+z_3\in\R$ and $z_1z_2z_3\in(-\infty,0)$. Then
\begin{equation}\label{eq:SU-three-angles}
 z_1+z_2+z_3\le \sum_{j=1}^3 |z_j|-2\min_j|z_j|.
\end{equation}
\end{lemma}
\begin{proof}
Write $z_j=a_je^{i\theta_j}$ with $a_j>0$ and
$\theta_j\in[-\pi,\pi]$. The hypotheses imply
\[
 \sum_j a_j\sin\theta_j=0,\qquad
 \theta_1+\theta_2+\theta_3\equiv\pi\pmod{2\pi}.
\]
We first show that $\sum_j\cos\theta_j\le1$.
If all the sines vanish, the three unit complex numbers are signs and
an odd number are negative; the assertion follows. Otherwise, the
vanishing weighted sine sum forces both positive and negative sines,
hence both positive and negative angles. Their sum, an odd multiple
of $\pi$, consequently equals $\pi$ or $-\pi$.
If the sum is $\pi$, precisely one angle is negative and two are
positive, and the trigonometric identity
\[
 \cos\theta_1+\cos\theta_2+\cos\theta_3-1
 =4\sin\frac{\theta_1}{2}\sin\frac{\theta_2}{2}
      \sin\frac{\theta_3}{2}
\]
shows the desired inequality. For a sum of $-\pi$, change the signs
of the angles. Thus $\sum_j\cos\theta_j\le1$ in all cases, whence
\[
 \begin{split}
  \sum_j\bigl(|z_j|-\operatorname{Re}z_j\bigr)
    &=\sum_j a_j(1-\cos\theta_j)\\
    &\ge \min_j a_j\sum_j(1-\cos\theta_j)
     \ge 2\min_j a_j.
 \end{split}
\]
The sum of the $z_j$ is real, which proves \eqref{eq:SU-three-angles}.
\end{proof}

\begin{lemma}\label{lem:SU-three-trace}
If $C\in\C_{3\times3}$ has singular values
$s_1\ge s_2\ge s_3>0$, $\Tr C\in\R$, and $\det C<0$,
then
\begin{equation}\label{eq:SU-three-trace}
 \Tr C\le s_1+s_2-s_3.
\end{equation}
\end{lemma}
\begin{proof}
Let $z_1,z_2,z_3$ be the eigenvalues of $C$, counted with
multiplicity. Lemma~\ref{lem:SU-three-complex} gives
\[
 \Tr C\le \sum_j|z_j|-2\min_j|z_j|.
\]
Schur triangularization and the trace-norm inequality yield
$\sum_j|z_j|\le\sum_j s_j$. Moreover, each eigenvalue $z_j$
satisfies $|z_j|\ge s_3$: for a normalized eigenvector $v$,
$|z_j|=\|Cv\|\ge s_3$. Combining these estimates proves
\eqref{eq:SU-three-trace}.
\end{proof}

\begin{theorem}[The special-unitary torus in order three]
\label{thm:SU-three}
Let $A\in\C_{3\times3}$ have singular values
$s_1\ge s_2\ge s_3\ge0$, and suppose its $\SU(3)\times \SU(3)$
orbit covers $\mathcal T(r)$. Then
\begin{equation}\label{eq:SU-three-bound}
 r_1+r_2+r_3\le s_1+s_2-s_3.
\end{equation}
Consequently \eqref{eq:SU-sufficient} is a necessary and sufficient
condition for full torus coverage in order three.
\end{theorem}
\begin{proof}
If $s_3=0$, the required inequality is the ordinary trace-norm
bound applied to any attainable diagonal. Suppose $s_3>0$ and write
$\det A=|\det A|e^{i\varphi}$. Choose $\theta\in\R$ with
$e^{i(\varphi-3\theta)}=-1$. Full torus coverage provides
$U,V\in \SU(3)$ with
\[
 \diag(UAV)=e^{i\theta}(r_1,r_2,r_3).
\]
For $C=e^{-i\theta}UAV$ we have
$\Tr C=r_1+r_2+r_3\in\R$ and
$\det C=-|\det A|<0$. Since its singular values are those of $A$,
Lemma~\ref{lem:SU-three-trace} proves \eqref{eq:SU-three-bound}.
The initial partial-sum conditions are necessary by Ky Fan; all
conditions in \eqref{eq:SU-sufficient} are sufficient by
Thompson's theorem \cite{Thompson1994}.
\end{proof}

\subsection{Failure of necessity in order four}

We first isolate a convenient reformulation for equal target moduli.

\begin{lemma}[Trace-circle criterion]\label{lem:SU-trace-circle}
Let $A\in\C_{n\times n}$ and $T\ge0$. The $\SU(n)\times \SU(n)$ orbit
of $A$ covers $\mathcal T((T/n,\ldots,T/n))$ if and only if
for every $\theta\in\R$ there exist $U,V\in \SU(n)$ with
\begin{equation}\label{eq:SU-trace-circle}
 \Tr(UAV)=Te^{i\theta}.
\end{equation}
\end{lemma}
\begin{proof}
The necessity follows by taking all the diagonal entries equal to
$(T/n)e^{i\theta}$. Conversely, assume \eqref{eq:SU-trace-circle}.
Every complex matrix $B$ is unitarily similar to a matrix with
constant diagonal $\Tr B/n$: this scalar lies in the numerical
range of $B$, and induction on the orthogonal compression supplies
an orthonormal basis with that constant diagonal.
The conjugating unitary may be multiplied by a scalar to have
determinant one, without changing the conjugation. Therefore, for
every $\theta$ the $\SU(n)\times \SU(n)$ orbit of $A$ contains a
matrix with diagonal $(T/n)e^{i\theta}(1,\ldots,1)$.
Given arbitrary $\lambda_j\in\mathbb T$, choose $\theta$ such that
$e^{in\theta}=\prod_j\lambda_j$. Left multiplication by
$\diag(\lambda_1e^{-i\theta},\ldots,\lambda_ne^{-i\theta})\in \SU(n)$
gives the desired diagonal vector.
\end{proof}

\begin{proposition}[The first failure occurs in order four]
\label{prop:SU-four}
Let
\[
 A=\diag(100,1,1,-1),\qquad
 r=(25.3,25.3,25.3,25.3).
\]
The $\SU(4)\times \SU(4)$ orbit of $A$ covers $\mathcal T(r)$, although
\begin{equation}\label{eq:SU-four-violation}
 \sum_{j=1}^4r_j=101.2>100+1+1-1=101.
\end{equation}
Consequently, Thompson's conjectured necessity of the terminal inequality
in \eqref{eq:SU-sufficient} holds for matrices of orders one, two, and
three, but fails in order four. Thus order four is the smallest dimension
admitting a counterexample.
\end{proposition}
\begin{proof}
Put $T=101.2$. By Lemma~\ref{lem:SU-trace-circle}, it suffices
to construct for every $\theta\in\R$ a matrix in the
$\SU(4)\times \SU(4)$ orbit of $A$ with trace $Te^{i\theta}$.
Choose $\delta\in[-\pi/3,\pi/3]$ so that
$3\delta\equiv\pi-4\theta\pmod{2\pi}$. For $|x|\le0.04$, put
\[
 q(x)=\frac{T-100\cos x}{\cos(\delta-x/3)},\qquad
 F(x)=100\sin x+q(x)\sin(\delta-x/3).
\]
The denominator in $q(x)$ is greater than $0.48$ and its numerator
lies in $[1.2,1.28)$. Hence $0<q(x)<1.28/0.48<3$ throughout the
interval. Since $100\sin(0.04)>3$, it follows that
$F(-0.04)<0<F(0.04)$. Choose $x$ with $F(x)=0$ and set
\[
 \alpha=\theta+x,\qquad \beta=\theta+\delta-x/3.
\]
There is $t\in\R$ satisfying $1+2\cos t=q(x)$. The diagonal matrix
\[
 B_\theta=\diag\bigl(100e^{i\alpha},
         e^{i(\beta+t)},e^{i(\beta-t)},e^{i\beta}\bigr)
\]
has singular values $(100,1,1,1)$ and determinant $-100$, since
$\alpha+3\beta=4\theta+3\delta\equiv\pi\pmod{2\pi}$.
The real and imaginary parts of $e^{-i\theta}\Tr B_\theta$
are $T$ and $0$ by the definitions of $q$ and $F$; hence
$\Tr B_\theta=Te^{i\theta}$.
Finally, $B_\theta=AV_\theta$ for a diagonal matrix
$V_\theta\in \SU(4)$, because $A$ and $B_\theta$ have the same
nonzero determinant. Lemma~\ref{lem:SU-trace-circle} establishes full torus coverage in order four.

The minimality of this dimension follows from the preceding results.
Thompson proved necessity in order two \cite[Section~3]{Thompson1994},
and Theorem~\ref{thm:SU-three} proves necessity in order three.
In order one, $\SU(1)$ is trivial, so full coverage of a circle
is possible only when $A=0$, for which the terminal inequality holds.
Thus order four is the first possible failure.
\end{proof}

\section*{AI disclosure}
During the preparation of this work, the authors used ChatGPT and Gemini to assist with language editing, exposition, and brainstorm. The authors are responsible for verifying the mathematical arguments and the final manuscript.

\section*{Funding}
Y. Luan is supported by the Shandong Provincial Natural Science Foundation (Grant No. ZR2025MS36).


\begin{thebibliography}{8}

\bibitem{BerensteinSjamaar2000}
Arkady Berenstein and Reyer Sjamaar, \emph{Coadjoint orbits, moment polytopes,
and the Hilbert--Mumford criterion}, Journal of the American Mathematical
Society \textbf{13} (2000), no.~2, 433--466.

\bibitem{GuoFan2022}
Cheng Guo and Shanhui Fan, \emph{Reciprocity constraints on reflection},
Physical Review Letters \textbf{128} (2022), 256101.

\bibitem{OSheaSjamaar2000}
Luis O'Shea and Reyer Sjamaar, \emph{Moment maps and Riemannian symmetric
pairs}, Mathematische Annalen \textbf{317} (2000), no.~3, 415--457.

\bibitem{Ressayre2010}
Nicolas Ressayre, \emph{Geometric invariant theory and the generalized
eigenvalue problem}, Inventiones Mathematicae \textbf{180} (2010),
no.~2, 389--441.

\bibitem{Tam1999}
Tin-Yau Tam, \emph{A Lie theoretic approach to Thompson's theorems on singular
values--diagonal elements and some related results}, Journal of the London
Mathematical Society (2) \textbf{60} (1999), no.~2, 431--448.

\bibitem{Thompson1977}
Robert~C. Thompson, \emph{Singular values, diagonal elements, and convexity},
SIAM Journal on Applied Mathematics \textbf{32} (1977), no.~1, 39--63.

\bibitem{Thompson1979}
Robert~C. Thompson, \emph{Singular values and diagonal elements of complex
symmetric matrices}, Linear Algebra and its Applications \textbf{26} (1979),
65--106.

\bibitem{Thompson1994}
Robert~C. Thompson, \emph{The diagonal torus of a matrix under special
unitary equivalence}, SIAM Journal on Matrix Analysis and Applications
\textbf{15} (1994), no.~3, 968--973.

\bibitem{Wootters1998}
William K. Wootters, \emph{Entanglement of formation of an arbitrary state
of two qubits}, Physical Review Letters \textbf{80} (1998), 2245--2248.

\end{thebibliography}
\end{document}